\documentclass[11pt,a4paper]{article}
\usepackage{amsmath}
\usepackage{amsfonts}
\usepackage{amssymb}
\usepackage{amsthm}
\usepackage{amscd}
\usepackage{graphicx}
\theoremstyle{plain}
\newtheorem{lemma}{Lemma}[section]
\newtheorem{theorem}{Theorem}[section]

\title{Lagrangian submanifolds with \\ prescribed second fundamental form}

\author{
	\emph{Bang-Yen Chen} \\
	Michigan State University \\
	Department of Mathematics \\ 
	East Lansing, MI 48824--1027, USA \\ 
	e-mail: \emph{bychen@math.msu.edu} \\
    \\ 
    \emph{Joeri Van der Veken}
	\thanks{Part-time postdoctoral fellow of the Research Foundation Flanders (F.W.O.)} \\
	KU Leuven \\
	Department of Mathematics \\ 
	Celestijnenlaan 200B -- Box 2400, BE--3001 Leuven, Belgium \\ 
	e-mail: \emph{joeri.vanderveken@wis.kuleuven.be} \\
    \\
	\emph{Luc Vrancken} \\
	Universit\'e de Valenciennes \\
	LAMAV, ISTV2 \\ 
	\smallskip
	Campus de Mont Houy, 59313 Valenciennes Cedex 9, France \\ 
	KU Leuven \\
	Department of Mathematics \\ 
	Celestijnenlaan 200B -- Box 2400, BE--3001 Leuven, Belgium \\ 
	e-mail: \emph{luc.vrancken@univ-valenciennes.fr} \\
	}

\date{}

\begin{document}

\maketitle

\begin{center} \emph{Dedicated to Franki Dillen, a great friend and colleague} \end{center}

\begin{abstract}
We classify Lagrangian submanifolds of complex space forms, whose second fundamental form can be written in a certain way, depending on a real parameter. For some special values of this parameter, the resulting submanifolds are ideal in the sense that they realize equality in an inequality for a delta-curvature.
\end{abstract}

\section{Introduction}

Let $\tilde M^n$ be a complex $n$-dimensional K\"ahler manifold with complex structure $J$. An isometric immersion $M^n \to \tilde M^n$ of a real $n$-dimensional Riemannian manifold $M^n$ into $\tilde M^n$ is called {\it Lagrangian\/} if $J$ is an isomorphism between the tangent and the normal space of the submanifold at any of its points. For instance the systems of partial differential equations of Hamilton-Jacobi type lead to the study of Lagrangian submanifolds and foliations in the cotangent bundle. Lagrangian submanifolds also play some important roles in supersymmetric field theories and in string theory. 

From now on, we will assume that the ambient space is a complex space form, i.e., a K\"ahler manifold of constant holomorphic sectional curvature $4c$, denoted by $\tilde M^n (4c)$. For $c=0$ this is the flat complex space $\mathbb C^n$, for $c>0$ it is a complex projective space $\mathbb C P^n(4c)$ and for $c<0$ it is a complex hyperbolic space $\mathbb C H^n(4c)$. 

We shall study Lagrangian submanifolds of complex space forms for which the second fundamental form takes a special form. In particular, we assume that every point has a neighbourhood on which there exists an orthonormal frame $\{e_1,\ldots,e_n\}$ and a function $\lambda$ such that the second fundamental form $h$ is given by
\begin{equation} \label{CVV:sff}
\begin{aligned}
& h(e_1,e_1) = \lambda Je_1, \\
& h(e_1,e_{\alpha}) = d \lambda Je_{\alpha}, \\
& h(e_{\alpha},e_{\beta}) = \delta_{\alpha\beta} d \lambda Je_1 + \sum_{\gamma=2}^n h_{\alpha\beta}^{\gamma} Je_{\gamma}, \mbox{ with } \sum_{\alpha = 2}^n h_{\alpha\alpha}^{\gamma} = 0 \mbox{ for every } \gamma = 2,\ldots,n.
\end{aligned}
\end{equation}
Here, $\alpha,\beta \in \{2,\ldots,n\}$ and $d$ is a real constant. There are two special values of $d$, for which the second fundamental form above appears in classifications of \emph{ideal} Lagrangian submanifolds. To explain what ideal Lagrangian submanifolds are, we need the notion of delta-invariants.

For the purpose of this paper, we just want to give a rough idea about delta-invariants and their inequalities. For a detailed overview, we refer to the book \cite{CVV:C}. On an $n$-dimensional Riemannian manifold $M^n$, one can define a delta-invariant $\delta(n_1,\ldots,n_k)$ for every $k$-tuple $(n_1,\ldots,n_k)$, satisfying the following properties: $2 \leq n_1 \leq \ldots \leq n_k \leq n-1$ and $n_1 + \ldots + n_k \leq n$. Any such delta-invariant is a curvature function on $M^n$. Now consider an immersion of $M^n$ into some ambient space, in particular, consider a Lagrangian immersion of $M^n$ into a complex space form $\tilde M^n(4c)$, as this is the case we are interested in. Then, for any $k$-tuple $(n_1,\ldots,n_k)$ as above, one can prove the following pointwise inequality:
$$ \delta(n_1,\ldots,n_k) \leq a(n,n_1,\ldots,n_k) \| H \|^2 + \frac 12 \left(n(n-1)-\sum_{j=1}^k n_j (n_j-1)\right) c, $$
where 
\begin{align*}
& a(n,n_1,\ldots,n_k) = \frac{n^2 \left( n-\sum_{j=1}^k n_j + 3k - 1 - 6 \sum_{j=1}^k \frac{1}{2+n_j} \right)}{2 \left(n -\sum_{j=1}^k n_j + 3k + 2 - 6 \sum_{j=1}^k \frac{1}{2+n_j} \right)} & \mbox{ if } n_1 + \ldots + n_k < n, \\
& a(n,n_1,\ldots,n_k) = \frac{n^2 \left( k-1-2 \sum_{j=2}^k \frac{1}{2+n_j} \right)}{2 \left( k-2 \sum_{j=2}^k \frac{1}{2+n_j} \right)} & \mbox{ if } n_1 + \ldots + n_k = n,
\end{align*}
cfr. \cite{CVV:CD}, \cite{CVV:CDVV1}. Moreover, the inequalities are optimal, in the sense that for every $\delta(n_1,\ldots,n_k)$, there exist non-minimal Lagrangian submanifolds realizing equality at every point. These submanifolds are called $\delta(n_1,\ldots,n_k)$-ideal. The two special choices for $d$ are now as follows.

\medskip

\textbf{Special case I:} {\mathversion{bold}$d = \frac{1}{2+m}$} \textbf{for some integer} {\mathversion{bold}$m \geq 2$} \textbf{dividing} {\mathversion{bold}$n-1.$} In this case, Lagrangian submanifolds with second fundamental form \eqref{CVV:sff} are $\delta(m,\ldots,m)$-ideal ($k$ times $m$), with $km = n-1$. The very particular case of $d = \frac{1}{n+1}$, corresponding to a $\delta(n-1)$-ideal Lagrangian submanifold, is treated in the paper \cite{CVV:CDV}.

\medskip

\textbf{Special case II:} {\mathversion{bold}$d = \frac{1}{n-1}$} \textbf{and} {\mathversion{bold}$n \geq 5.$} In this case, some Lagrangian submanifolds with second fundamental form \eqref{CVV:sff} are $\delta(2,n-2)$-ideal. In fact, an extra condition on $h$ involving the existence of a second special direction, orthogonal to $e_1$, is needed. The resulting second fundamental form describes only one of the many possible cases for the second fundamenal form of a $\delta(2,n-2)$-ideal submanifold, as is discussed in the paper \cite{CVV:CDVV2}. 

In the present paper, we give a classification result for Lagrangian submanifolds of complex space forms with second fundamental form \eqref{CVV:sff} for arbitrary values of the constant $d \neq 0, \frac 12$. If $\lambda = 0$, the submanifold is minimal and $\delta(n-1)$-ideal and thus we will assume that $\lambda \neq 0$. We obtain that $M^n$ is intrinsically a warped product $I \times_f N$ of an interval $I$ and an $(n-1)$-dimensional manifold $N$. Here, $e_1$ is tangent to $I$, whereas $e_2,\ldots, e_n$ are tangent to $N$. The immersion is then built out of two things. On one hand a curve depending on a parameter $t \in I$, which is determined by a system of ODE's, and on the other hand a minimal Lagrangian immersion of $N$ for which the components of the second fundamental form are, up to a factor depending on $t$, equal to the corresponding components of $h$. Some fundamental lemmas leading to this classification result are given in Section 2, whereas explicit expressions for the immersions are given in Section 3. The three theorems in this section correspond to the cases that the ambient space is $\mathbb C^n$, $\mathbb C P^n(4)$ or $\mathbb C H^n(-4)$ respectively. We ende the paper with some remarks on special cases.

\section{Some fundamental lemmas}

Consider a Lagrangian submanifold $M^n$ of a complex space form of constant holomorphic sectional curvature $c$, with second fundamental form $h$. We denote by $\nabla$ the Levi-Civita connection of the submanifold. The normal connection is then completely determined by $\nabla^{\perp}_X JY = J \nabla_X Y$ for vector fields $X$, $Y$ tangent to $M^n$. A key property of Lagrangian submanifolds is that the cubic form $\langle h(X,Y) , JZ \rangle$ is totally symmetric in the tangent vectors $X,Y,Z$. The fundamental equations of Gauss and Codazzi reduce respectively to
\begin{equation} \label{CVV:gauss}
\langle R(X,Y)Z,W \rangle \!=\! \langle h(X,W),h(Y,Z) \rangle \!-\! \langle h(X,Z),h(Y,W) \rangle 
+ c(\langle X,W \rangle \langle Y,Z \rangle \!-\! \langle X,Z \rangle \langle Y,W \rangle)
\end{equation}
and
\begin{equation} \label{CVV:codazzi}
(\nabla h)(X,Y,Z) = (\nabla h)(Y,X,Z),
\end{equation}
where $R$ is the Riemann-Christoffel curvature tensor of $M^n$ and the covariant derivative of $h$ is defined by $(\nabla h)(X,Y,Z) = \nabla^{\perp}_X h(Y,Z) - h(\nabla_X Y,Z) - h(Y,\nabla_X Z)$.

As mentioned above, if the second fundemantal form of a Lagrangian submanifold of a complex space form is given by \eqref{CVV:sff} with $\lambda = 0$, the submanifold is minimal and $\delta(n-1)$-ideal. Hence, we will assume from now on that $\lambda \neq 0$. In order to get a classification, we additionally assume that $d \neq 0, \frac 12$. Remark that this is satisfied for the special cases described in the introduction. The proofs of the lemmas below are similar to those given in \cite{CVV:CDV}, so we keep them rather short.

\begin{lemma} \label{CVV:lem1}
Let $M^n \to \tilde M^n(4c)$ be a Lagrangian submanifold of a complex space form of dimension $n \geq 3$. Assume that there is a local orthonormal frame $\{e_1,\ldots,e_n\}$ on $M^n$ with respect to which the second fundamental form of the immersion takes the form \eqref{CVV:sff} with $\lambda \neq 0$ and $d \neq 0, \frac 12$. Then
\begin{itemize}
\item[(i)] $\nabla_{e_1}e_1 = 0$,
\item[(ii)] $e_{\alpha}(\lambda) = 0$ for $\alpha = 2,\ldots,n$,
\item[(iii)] $\nabla_{e_{\alpha}}e_1 = \mu e_{\alpha}$, with $\mu = \frac{d}{(1-2d)\lambda}e_1(\lambda)$ for $\alpha = 2,\ldots,n$.
\end{itemize}
\end{lemma}

\emph{Proof.} Assume that $\nabla_{e_1}e_1 \neq 0$. Then we can choose $e_2$ such that $\nabla_{e_1}e_1 = b e_2$ for some $b \neq 0$. We will use the equation of Codazzi several times. From $(\nabla h)(e_{\alpha},e_1,e_1) = (\nabla h)(e_1,e_{\alpha},e_1)$, where $\alpha \geq 2$, one finds by a straightforward computation
\begin{equation*}
e_2(\lambda) = (1-2d) b \lambda, \quad
e_{\beta}(\lambda) = 0 \mbox{ for $\beta \geq 3$, } \quad
(1-2d) \lambda \nabla_{e_{\alpha}}e_1 = d e_1(\lambda) e_{\alpha} + b J h(e_2,e_{\alpha})^{\perp},
\end{equation*}
where $h(X,Y)^{\perp}$ denotes the component of $h(X,Y)$ perpendicular to $Je_1$. Using these, one finds from
$$ \sum_{\alpha=3}^n \langle (\nabla h)(e_2,e_{\alpha},e_1) - (\nabla h)(e_{\alpha},e_2,e_1), Je_{\alpha} \rangle = 0 $$
that
\begin{equation} \label{CVV:eq1}
\sum_{\alpha=2}^n \langle h(e_2,e_{\alpha})^{\perp},h(e_2,e_{\alpha})^{\perp} \rangle = (n-2)d(1-2d)^2\lambda^2
\end{equation}
and from
$$ \sum_{\alpha=2}^n \langle (\nabla h)(e_1,e_{\alpha},e_{\alpha}) - (\nabla h)(e_{\alpha},e_1,e_{\alpha}) , Je_2 \rangle = 0 $$
that
\begin{equation} \label{CVV:eq2}
\sum_{\alpha=2}^n \langle h(e_2,e_{\alpha})^{\perp},h(e_2,e_{\alpha})^{\perp} \rangle = (n+2d)d(1-2d)\lambda^2
\end{equation}
By comparing \eqref{CVV:eq1} and \eqref{CVV:eq2} and using the fact that $d \neq 0, \frac 12$, we obtain that $d = \frac{1}{1-n}$. But this is impossible since \eqref{CVV:eq1} implies that $d$ is non-negative. Hence, we conclude that our first assumption was wrong and that $\nabla_{e_1}e_1 = 0$. The two other equations now follow from $(\nabla h)(e_{\alpha},e_1,e_1) = (\nabla h)(e_1,e_{\alpha},e_1)$ for $\alpha \geq 2$.
\hfill $\square$ 

\begin{lemma} \label{CVV:lem2}
With the notations and the assumptions of Lemma \ref{CVV:lem1}, we have
\begin{itemize}
\item[(i)] $\mathcal D_1 = \mathrm{span}\{e_1\}$ is a totally geodesic foliation,
\item[(ii)] $\mathcal D_2 = \mathrm{span}\{e_2,\ldots,e_n\}$ is a spherical foliation,
\item[(iii)] $M^n$ is locally a warped product $I \times_f N$ of an integral curve of $e_1$ and a leaf of $\mathcal D_2$.
\end{itemize}
\end{lemma}

\emph{Proof.} By a result of Hiepko, see for example \cite{CVV:C}, p. 90, it is sufficient to prove the first two statements of the lemma to conclude the last one.

The distribution $\mathcal D_1$ is integrable since it is one-dimensional. It follows from Lemma \ref{CVV:lem1} \emph{(i)} that the leaves of the resulting foliation are geodesics.

The integrability of $\mathcal D_2$ follows from Lemma \ref{CVV:lem1} \emph{(iii)}. It remains to prove that the resulting foliation is spherical, i.e., that the leaves are totally umbilical with parallel mean curvature vector. Total umbilicity follows immediately from Lemma \ref{CVV:lem1} \emph{(iii)}, and it only remains to verify that $e_{\alpha}(\mu)=0$ for $\alpha \geq 2$. For this purpose, consider the system
\begin{equation} \label{sys}
\left \{ \begin{array}{l} e_1(\lambda) =  \displaystyle{\frac{1-2d}{d}} \lambda\mu, \vspace{.2cm} \\ e_{\alpha}(\lambda)=0, \end{array} \right.
\end{equation}
which is obtained from Lemma \ref{CVV:lem1} \emph{(ii), (iii)}. It follows from Lemma \ref{CVV:lem1} \emph{(i), (iii)} that $[e_1,e_{\alpha}]$ is perpendicular to $e_1$ and hence $[e_1,e_{\alpha}](\lambda)=0$. The integrability condition for the system \eqref{sys} now reduces to $e_{\alpha}(\mu)=0$. \hfill $\square$ 

\bigskip

If we choose a local coordinate system $(t,u_2,\ldots,u_n)$ on $M^n = I \times_f N$ such that $e_1 = \frac{\partial}{\partial t}$, then we have the following.

\begin{lemma} \label{CVV:lem3}
With the notations and the assumptions of Lemma \ref{CVV:lem1}, the functions $\lambda$ and $\mu$ are functions of $t$ only and they satisfy the following system of ODE's:
\begin{align*}
& \lambda' = \frac{1-2d}{d} \lambda \mu,
& \mu' = -c -\mu^2 -d(1-d)\lambda^2.
\end{align*}
\end{lemma}

\emph{Proof.} The fact that $\lambda$ only depends on $t$ follows from Lemma \ref{CVV:lem1}, \emph{(ii)} and the fact that $\mu$ only depends on $t$ follows from the fact that $\mathcal D_2$ is spherical (cfr. proof of Lemma \ref{CVV:lem2}). The first differential equation is just a reformulation of the definition of $\mu$, given in Lemma \ref{CVV:lem1} \emph{(iii)}. The second equation follows from the equation of Gauss \eqref{CVV:gauss}. In particular, one can check that $R(e_{\alpha},e_1)e_1 = \nabla_{e_{\alpha}}\nabla_{e_1}e_1 - \nabla_{e_1}\nabla_{e_{\alpha}}e_1 -\nabla_{[e_{\alpha},e_1]}e_1 = - e_1(\mu)e_{\alpha} - \mu^2e_{\alpha}$, and hence the equation
$$ \langle R(e_{\alpha},e_1)e_1,e_{\alpha} \rangle = \langle h(e_1,e_1),h(e_{\alpha},e_{\alpha})\rangle - \langle h(e_1,e_{\alpha}),h(e_1,e_{\alpha}) \rangle + c $$
is equivalent to the second ODE in the lemma.
\hfill $\square$ 

\section{The classification results}

In this section, we will describe the Lagrangian immersions with second fundamental form \eqref{CVV:sff} as explicitly as possible. In the case of $\mathbb C P^n(4)$ and $\mathbb C H^n(-4)$, we will use \emph{horizontal lifts} to do this. 

Consider the unit sphere $S^{2n+1}(1) = \{ (z_1,z_2,\ldots,z_{n+1}) \in \mathbb C^{n+1} \ | \ |z_1|^2 + |z_2|^2 + \ldots + |z_{n+1}|^2 = 1 \}$ equipped with the standard Sasakian structure. In particular, the contact form $\eta$ is dual to the vector field defined by $\xi(z) = iz$ for all $z \in S^{2n+1}(1)$. A submanifold $N$ of $S^{2n+1}(1)$ is called Legendrian if $\eta(X)=0$ for all $X$ tangent to $N$. Now let $\pi : S^{2n+1}(1) \to \mathbb C P^n(4)$ be the Hopf fibration. Recall that the fibers of $\pi$ are precisely the integral curves of $\xi$. It is well-known that every Lagrangian immersion of a manifold $N$ into $\mathbb C P^n(4)$ can be lifted locally, or globally if $N$ is simply connected, to a Legendrian immersion of the same manifold into $S^{2n+1}(1)$. 

Lagrangian immersions into $\mathbb C H^n(-4)$ can be treated in a similar way. For this purpose, consider $H^{2n+1}_1(-1) = \{ (z_1,z_2,\ldots,z_{n+1}) \in \mathbb C^{n+1}_1 \ | \ -|z_1|^2 + |z_2|^2 + \ldots + |z_{n+1}|^2 = -1 \}$ with the induced Lorentzian metric. As in the previous case, the vector field defined by $\xi(z)=iz$ for every $z \in H^{2n+1}_1(-1)$ is tangent to $H^{2n+1}_1(-1)$ and satisfies $\langle \xi,\xi \rangle = -1$. The orthogonal complement of $\xi(z)$ in $T_zH^{2n+1}_1(-1)$ has a positive definite metric, and there exists a Riemannian submersion $\pi : H^{2n+1}_1(-1) \to \mathbb C H^n(-4)$, analogous to the Hopf fibration. Again, a submanifold $N$ of $H^{2n+1}_1(-1)$ will be called Legendrian if $\eta(X)=0$ for every $X$ tangent to $N$, where $\eta$ is the dual one-form of $\xi$. Again, every Lagrangian immersion of a manifold $N$ into $\mathbb C H^n(-4)$ can be lifted locally, or globally if $N$ is simply connected, to a Legendrian immersion of the same manifold into $H^{2n+1}_1(-1)$. 

For more details on horizontal lifts of Lagrangian submanifolds, we refer to \cite{CVV:R}

The proofs of our main theorems are very similar to those of the main results in \cite{CVV:CDV}, where a very special case of the first special case mentioned in the introduction was treated, and we will therefore omit them. 

In all three theorems, the function $\lambda$ is the function appearing in the second fundamental form \eqref{CVV:sff} and the function $\mu$ is the one introduced in Lemma \ref{CVV:lem1}. 

\begin{theorem} \label{CVV:theo1}
Let $M^n \to \mathbb C^n$ be a Lagrangian submanifold of dimension $n \geq 3$. Assume that there is a local orthonormal frame $\{e_1,\ldots,e_n\}$ on $M^n$ with respect to which the second fundamental form of the immersion takes the form \eqref{CVV:sff} with $\lambda \neq 0$ and $d \neq 0, \frac 12$. Then, up to rigid motions of $\mathbb C^n$, the immersion is given by
$$ L(t,u_2,\ldots,u_n) = \frac{e^{i\theta}}{\mu + id\lambda} \phi(u_2,\ldots,u_n). $$
Here, $\frac{\partial}{\partial t}=e_1$ and $\lambda$, $\mu$ and $\theta$ are functions of $t$ only, satisfying
$$ \lambda' = \frac{1-2d}{d} \lambda\mu, \qquad \mu' = -\mu^2 - d(1-d)\lambda^2, \qquad \theta' = \lambda. $$
Moreover, $\phi$ is a minimal Legendrian immersion into $S^{2n-1}(1) \subseteq \mathbb C^n$.
\end{theorem}

Remark that the system of differential equations in Theorem \ref{CVV:theo1} can be explicitly solved. The general solution depends on whether $d=1$ or $d\neq 1$. For $d=1$, the solution is given by
$$ \mu(t) = \frac{1}{t+k_1}, \qquad \lambda=\frac{k_2}{t+k_1}, \qquad \theta(t) = k_2 \ln |t+k_1| + k_3. $$
After a reparametrization and an isometry, one can assume that $k_1=k_3=0$.
Now assume that $d \neq 1$. After replacing $e_1$ by $-e_1$ if necessary, one can assume that $\lambda$ is positive. It then follows from the differential equations that $\lambda^{\frac{2d}{1-2d}}(\mu^2+d^2\lambda^2)$ is a positive constant, say $k^2$. Then $\mu = \pm\sqrt{k^2 \lambda^{\frac{2d}{2d-1}} - d^2 \lambda^2}$, and one obtains
$$ \frac{\mathrm d\theta}{\mathrm d\lambda} = \pm \frac{d}{1-2d} \frac{1}{\sqrt{k^2 \lambda^{\frac{2d}{2d-1}} - d^2 \lambda^2}}, $$
for which the solution is given by
$$ \theta(\lambda) = \pm \frac{1}{d-1} \, \mathrm{csc}^{-1}\left(\frac kd \lambda^{-\frac{d-1}{2d-1}}\right). $$
One does not have to write an additive constant, since this can be assumed to be zero after an isometry. After the substitution $t \mapsto \lambda(t)$, the expressions for $\mu(\lambda)$ and $\theta(\lambda)$ given above, determine the coefficient of $\phi$ in the immersion in Theorem \ref{CVV:theo1}.

\begin{theorem} \label{CVV:theo2}
Let $M^n \to \mathbb C P^n(4)$ be a Lagrangian submanifold of dimension $n \geq 3$. Assume that there is a local orthonormal frame $\{e_1,\ldots,e_n\}$ on $M^n$ with respect to which the second fundamental form of the immersion takes the form \eqref{CVV:sff} with $\lambda \neq 0$ and $d \neq 0, \frac 12$. Then, up to isometries of $\mathbb C P^n(4)$, the horizontal lift of the immersion to $S^{2n+1}(1) \subseteq \mathbb C^{n+1} = \mathbb C^n \times \mathbb C$ is given by
$$ L(t,u_2,\ldots,u_n) = \left( \frac{e^{id\theta} \phi(u_2,\ldots,u_n)}{\sqrt{1 + \mu^2 + d^2\lambda^2}}, \frac{e^{i(1-d)\theta} (id\lambda - \mu)}{\sqrt{1 + \mu^2 + d^2\lambda^2}} \right). $$
Here, $\frac{\partial}{\partial t}=e_1$ and $\lambda$, $\mu$ and $\theta$ are functions of $t$ only, satisfying
$$ \lambda' = \frac{1-2d}{d} \lambda\mu, \qquad \mu' = -1 -\mu^2 - d(1-d)\lambda^2, \qquad \theta' = \lambda. $$
Moreover, $\phi$ is a minimal Legendrian immersion into $S^{2n-1}(1) \subseteq \mathbb C^n$.
\end{theorem}

\begin{theorem} \label{CVV:theo3}
Let $M^n \to \mathbb C H^n(-4)$ be a Lagrangian submanifold of dimension $n \geq 3$. Assume that there is a local orthonormal frame $\{e_1,\ldots,e_n\}$ on $M^n$ with respect to which the second fundamental form of the immersion takes the form \eqref{CVV:sff} with $\lambda \neq 0$ and $d \neq 0, \frac 12$. Then, up to isometries of $\mathbb C H^n(-4)$, the horizontal lift of the immersion to $H^{2n+1}_1(-1) \subseteq \mathbb C^{n+1}_1$ depends on the sign of $1- \mu^2 - d^2 \lambda^2$. In each case, $\frac{\partial}{\partial t}=e_1$.
\begin{itemize}
\item If $\mu^2 + d^2 \lambda^2 < 1$, then
$$ L(t,u_2,\ldots,u_n) = \left( \frac{e^{id\theta} \phi(u_2,\ldots,u_n)}{\sqrt{1 - \mu^2 - d^2\lambda^2}}, \frac{e^{i(1-d)\theta} (id\lambda - \mu)}{\sqrt{1 - \mu^2 - d^2\lambda^2}} \right), $$
where $\lambda$, $\mu$ and $\theta$ are functions of $t$ only, satisfying
$$ \lambda' = \frac{1-2d}{d} \lambda\mu, \qquad \mu' = 1 -\mu^2 - d(1-d)\lambda^2, \qquad \theta' = \lambda $$
and $\phi$ is a minimal Legendrian immersion into $H^{2n-1}_1(-1) \subseteq \mathbb C^n_1$.
\item If $\mu^2 + d^2 \lambda^2 > 1$, then
$$ L(t,u_2,\ldots,u_n) = \left( \frac{e^{i(1-d)\theta} (id\lambda - \mu)}{\sqrt{\mu^2 + d^2\lambda^2 - 1}}, \frac{e^{id\theta} \phi(u_2,\ldots,u_n)}{\sqrt{\mu^2 + d^2\lambda^2 - 1}} \right), $$
where $\lambda$, $\mu$ and $\theta$ are functions of $t$ only, satisfying
$$ \lambda' = \frac{1-2d}{d} \lambda\mu, \qquad \mu' = 1 -\mu^2 - d(1-d)\lambda^2, \qquad \theta' = \lambda $$
and $\phi$ is a minimal Legendrian immersion into $S^{2n-1}(1) \subseteq \mathbb C^n$.
\item If $\mu^2 + d^2 \lambda^2 = 1$, then
\begin{multline*}
L(t,u_2,\ldots,u_n) = \frac{e^{\frac{2d}{2d-1} i \tan^{-1}\left((\tanh(\frac{1-2d}{2d}t)\right)}}{\cosh^{\frac{2d}{2d-1}}\left(\frac{1-2d}{2d}t\right)} \left[ \left(w + \frac i2 \langle \phi,\phi \rangle + i, \phi, w + \frac i2 \langle \phi,\phi \rangle \right) \right.\\
+ \left. \int_0^t \cosh^{\frac{2d}{2d-1}} \left(\frac{1-2d}{2d}t \right)e^{2i \tan^{-1}\left(\tanh(\frac{1-2d}{2d}t)\right)}dt \ (1,0,\ldots,0,1)\right],
\end{multline*}
where $\phi$ is a minimal Lagrangian immersion into $\mathbb C^{n-1}$ and $w$ is a solution of the system
$$ \frac{\partial w}{\partial u_{\alpha}} = \langle i \frac{\partial \phi}{\partial u_{\alpha}}, \phi \rangle, \mbox{ for } \alpha=2,\ldots,n. $$
\end{itemize}
\end{theorem}

\section{Some final remarks}

Throughout the paper we assumed that $n \geq 3$. If $n=2$, one can check, by similar computations, that the conclusions of Lemma \ref{CVV:lem1} are still true, except if $d=-1$. The case $d=-1$ corresponds to a minimal Lagrangian surface since the second fundamental form \eqref{CVV:sff} reduces to $h(e_1,e_1)=\lambda Je_1$, $h(e_1,e_2)=-\lambda Je_2$ and $h(e_2,e_2)=-\lambda Je_1$. In the case $n=2$, $d \neq -1$, all the following lemmas and theorems are still true. Moreover, the Legendrian immersions $\phi$ can be explicitly given: in Theorem \ref{CVV:theo1} and Theorem \ref{CVV:theo2}, $\phi$ is a horizontal geodesic in $S^3(1) \subseteq \mathbb C^2$ and thus, after an isometry, $\phi(u_2) = (\cos u_2, \sin u_2)$; in Theorem \ref{CVV:theo3} one has, by a similar argument respectively $\phi(u_2) = (\cosh u_2, \sinh u_2)$, $\phi(u_2) = (\cos u_2,\sin u_2)$, $\phi(u_2)=u_2$ and $w$ is constant.

In the other excluded cases, namely $d=0$ and $d = \frac 12$, it is not clear that $\nabla_{e_1}e_1 = 0$ and that the submanifold will have a warped product structure.

%

\end{document}